\documentclass[doc]{apa6}

\usepackage[american]{babel}
\usepackage{csquotes}
\usepackage{cite}

\usepackage{graphicx}
\usepackage[T1]{fontenc}
\usepackage[utf8x]{inputenc}
\usepackage{calc}
\usepackage{graphicx}
\usepackage[unicode=true]
 {hyperref}
\usepackage{breakurl}
\usepackage{xcolor}
\usepackage[normalem]{ulem}

\usepackage{skak}
\usepackage{chessfss}

\usepackage{amsmath,amssymb}

\author{Zolt\'an Kov\'acs}

\affiliation{Private P\"adagogische Hochschule der Di\"ozese Linz, A-4020 Linz, Salesianumweg 3}

\title{The epsilon-delta game and the chess}
\shorttitle{The epsilon-delta game and the chess}

\abstract{The limit of a sequence by the
definition with $\varepsilon$ is introduced by the notion of checkmate in two moves.
The idea is also extended to
define the limit of a function with $\varepsilon$ and $\delta$.}


\begin{document}
\maketitle

\section{Introduction}
{\it Understanding} the definition of the limit of sequences and functions
remains a challenge across generations. For the teacher it is also an effort
and an intellectual challenge to introduce the essence
and content of the definition, including choice of examples for
their audience, enabling understanding.
Students are different, of course, but by broadening of the education
the delivery of mathematical content for the widest possible audience demands effective teaching
methods.

In this paper an idea will be suggested to define the limit of a sequence by the
definition with $\varepsilon$. This idea is not really new, what is more,
it may be a typical way to define it as a game of two players. In this article
this game will be more precisely described. This method can also be extended to
define the limit of a \textit{function} with $\varepsilon$ and $\delta$.

Interactive materials are also available on the Internet
for the epsilon-delta game. One of the first published materials is \cite{Townsend2000},
a Mathematica applet, was announced 15 years ago. Other example is
the GeoGebra material \cite{Mulholland2013} which can be directly launched
in a web browser.

\section{The limit of a sequence}

One typical definition for the notion of limit is as follows: 
Let the real number sequence $a_n$ be given. We call the real number $a$ the
limit of the sequence $a_n$ if for each positive number $\varepsilon$ there exists a natural number $N$ such that
for $n>N$ we have $|a-a_n|<\varepsilon$. Or, symbolically:
$$\exists a\in {\bf R}\hskip0.7cm
\forall \varepsilon >0\hskip0.7cm
\exists N\in {\bf N}\hskip0.7cm
\forall n>N (n\in {\bf N})\hskip1cm
|a-a_n|<\varepsilon.
$$
The first existential quantifier is not always highlighted, but in the following
it will have a key role.

When teaching this definition, for the very first approach the epsilon will have
concrete values, for which concrete numbers for $N$ will be searched. This will
make the definition shorter:
$$
\exists N\in {\bf N}\hskip0.7cm
\forall n>N (n\in {\bf N})\hskip1cm
|a-a_n|<\varepsilon.
$$

While the general definition can be described with the scheme
$$\exists\hskip0.3cm
\forall\hskip0.3cm
\exists\hskip0.3cm
\forall$$
here the scheme
$$\exists\hskip0.3cm
\forall
$$
can be found. Actually these schemes are identical with the structure
of the final one or two moves of a game of two
players, Alice and Bob,
\begin{enumerate}
\item Alice has a move such that
\item for any moves of Bob's,
\item Alice has another move such that,
\item for any other moves of Bob's, he loses;
\end{enumerate}
moreover
\begin{enumerate}
\item Alice has a move such that,
\item for any other moves of Bob's, he loses.
\end{enumerate}

From the perspective of mathematical logic actually any type of game of two players could
be acceptable to demonstrate the structure as an every day example, there are educational benefits
to prefer using the chess game here:
\begin{itemize}
\item Chess is an essential element of human culture given a chessboard is recognizable by almost everyone (even if the rules of chess are not fully understood).
\item The rules of chess are mostly understood by students. And, where this is not the case it is possible to introduce the general rules of the game through a short lesson which will be enough to understand the definition of limit in mathematics (see below).
\item The game is visual, two dimensional, hence the focus will be concentrated on a different interpretation: this is {\it graphic} representation of knowledge to extend the {\it numeric} and {\it verbal} reasoning \cite{Deanin2003}.
\item In many countries there are deep traditions in teaching chess in schools.
\item By the analysis of easier or even more difficult chess puzzles the structure of the definition of the limit
can be simplified to observe just a finite number of cases.
\end{itemize}

\begin{figure}
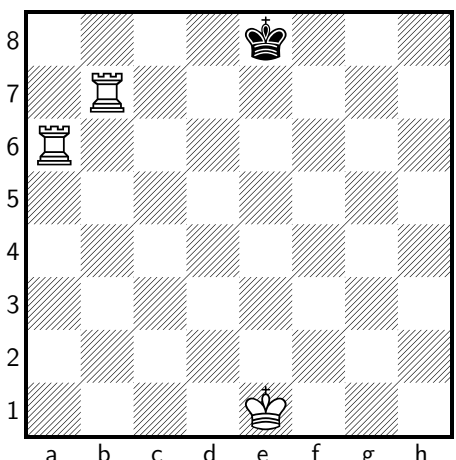

\centering
\fenboard{4k3/1R6/R7/8/8/8/8/4K3 w - - 0 20}
\showboard
\caption{Mate in one}
\label{fig:mate-in-one}
\end{figure}

\begin{figure}
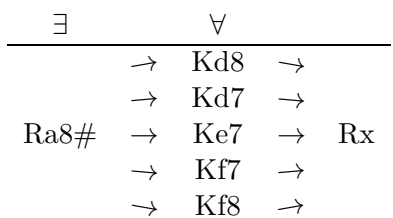

\centering
\begin{tabular}{ccccc}
$\exists$ &                                         & $\forall$ &               &    \\
\hline
          & \rotatebox[origin=c]{10}{$\rightarrow$}  & Kd8       & \rotatebox[origin=c]{-10}{$\rightarrow$} & \\
          & \rotatebox[origin=c]{5}{$\rightarrow$}   & Kd7       & \rotatebox[origin=c]{-5}{$\rightarrow$} & \\
Ra8\#     & \rotatebox[origin=c]{0}{$\rightarrow$}   & Ke7       & $\rightarrow$ & Rx \\
          & \rotatebox[origin=c]{-5}{$\rightarrow$}  & Kf7       & \rotatebox[origin=c]{5}{$\rightarrow$} &\\
          & \rotatebox[origin=c]{-10}{$\rightarrow$} & Kf8       & \rotatebox[origin=c]{10}{$\rightarrow$} & \\
\end{tabular}
\caption{Mate in one, example graph of possible moves}
\label{fig:mate-in-one-graph}
\end{figure}

Fig.~\ref{fig:mate-in-one} demonstrates the idea described above.
It is possible to make a connection between finding $N$ for a given epsilon and solving
a ``mate in one'' chess puzzle. For those students who do not know the rules of chess, a simple
puzzle can be shown. On the chessboard (see the figure above) in this example there are only two kinds of figures:
kings and rooks with the moves easily explained. In our case for a concrete $\varepsilon$
the appropriate $N$ is the rook move a8 (Ra8), so that the black king cannot have a proper move
anymore: for all 5 ``possible'' moves the black king remains in check, hence he is in checkmate.
This means that the move Ra8 corresponds to the definition's $N$, so that any move of Black
results in ``capturing the black king'', that is, checkmate (see Fig.~\ref{fig:mate-in-one-graph}).

For those students who already know the rules of chess, it is possible to solve more difficult puzzles.
There are plenty of examples available on the Internet about such puzzles, one possible site is \url{mateinone.com}.

\begin{figure}
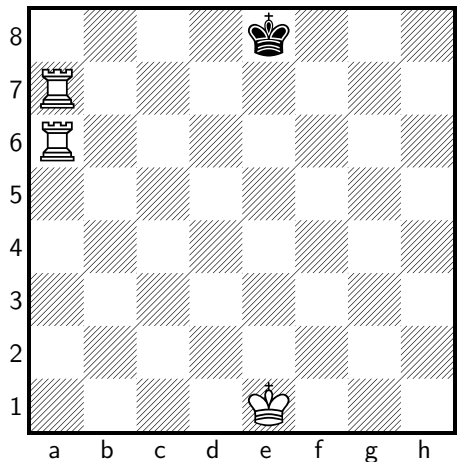

\centering
\fenboard{4k3/R7/R7/8/8/8/8/4K3 w - - 0 20}
\showboard
\caption{Mate in two}
\label{fig:mate-in-two}
\end{figure}

\begin{figure}
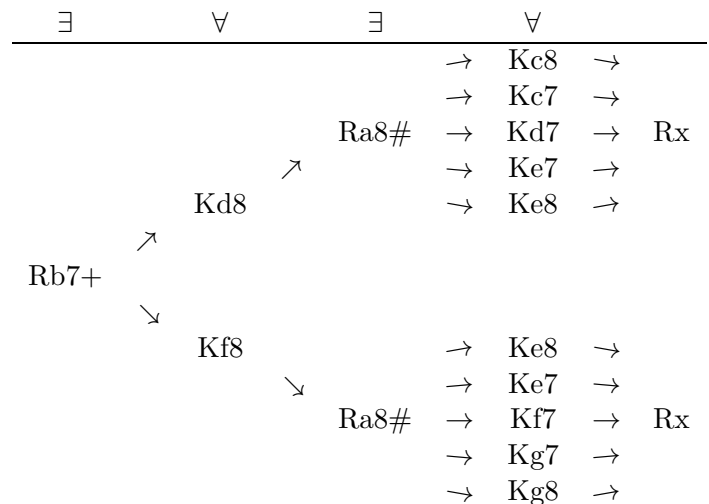

\centering
\begin{tabular}{ccccccccc}
$\exists$ &                                         & $\forall$ &               & $\exists$ &                                       & $\forall$ &               &    \\
\hline
     &                                          &          &                                                      &           & \rotatebox[origin=c]{10}{$\rightarrow$}  & Kc8       & \rotatebox[origin=c]{-10}{$\rightarrow$} & \\
     &                                          &          &                                                      &           & \rotatebox[origin=c]{5}{$\rightarrow$}   & Kc7       & \rotatebox[origin=c]{-5}{$\rightarrow$} & \\
     &                                          &          &                                                      & Ra8\#     & \rotatebox[origin=c]{0}{$\rightarrow$}   & Kd7       & $\rightarrow$ & Rx \\
     &                                          &          & \rotatebox[origin=c]{45}{$\rightarrow$}              &           & \rotatebox[origin=c]{-5}{$\rightarrow$}  & Ke7       & \rotatebox[origin=c]{5}{$\rightarrow$} &\\
     &                                          & Kd8      &                                                      &           & \rotatebox[origin=c]{-10}{$\rightarrow$} & Ke8       & \rotatebox[origin=c]{10}{$\rightarrow$} & \\
     &  \rotatebox[origin=c]{45}{$\rightarrow$} &          &                                                      &           &                                          &           &                                         & \\
Rb7+ &                                          &          &                                                      &           &                                          &           &                                         & \\
     &  \rotatebox[origin=c]{-45}{$\rightarrow$}&          &                                                      &           &                                          &           &                                         & \\
     &                                          & Kf8      &                                                      &           & \rotatebox[origin=c]{10}{$\rightarrow$}  & Ke8        &\rotatebox[origin=c]{-10}{$\rightarrow$}  & \\
     &                                          &          & \rotatebox[origin=c]{-45}{$\rightarrow$}             &           & \rotatebox[origin=c]{5}{$\rightarrow$}  & Ke7        & \rotatebox[origin=c]{-5}{$\rightarrow$}& \\
     &                                          &          &                                                      & Ra8\#     & \rotatebox[origin=c]{0}{$\rightarrow$}   & Kf7       & \rotatebox[origin=c]{-0}{$\rightarrow$}& Rx\\
     &                                          &          &                                                      &           & \rotatebox[origin=c]{-5}{$\rightarrow$}   & Kg7       & \rotatebox[origin=c]{5}{$\rightarrow$} &\\
     &                                          &          &                                                      &           & \rotatebox[origin=c]{-10}{$\rightarrow$}  & Kg8       & \rotatebox[origin=c]{10}{$\rightarrow$} & \\
\end{tabular}
\caption{Mate in two, example graph of possible moves}
\label{fig:mate-in-two-graph}
\end{figure}

Our students may need more time to prescind from the ``mate in one'' to the ``mate in two''. It can
be a great intellectual challenge for many students to change the concrete $\varepsilon$ to a general
one (see for example \cite{SwinyardLarsen2012}). A minor change of the figures on our chessboard (see the figure below)
will definitely change the puzzle, since White cannot deliver checkmate anymore in one move, that player needs
at least two moves (see Fig.~\ref{fig:mate-in-two}).
Indeed, one of the rooks must be moved to the vertical column b, and after any response of the black king
the rook on horizontal row 6 must be moved to row 8. The question may be ``reduced
to a previously solved problem'' by suggesting the moves Rb7--K?--Ra8.
Finally the moves described as Rb7--K?--Ra8--K? are the same as the structure $\exists\forall\exists\forall$
which will prepare the students for the precise definition of the limit (see Fig.~\ref{fig:mate-in-two-graph}).

It is clear that for White's moves there are multiple arrows to continue the game (since {\it for all} moves of Black's there exists a response of White's),
while Black's moves are continued only by one arrow (since it is sufficient for White to find {\it one} suitable move).

More difficult puzzles can be found at e.g.~\url{www.chesspuzzles.com/mate-in-two}.

\section{The limit of a function}

Obviously, the definition of the limit of a function has the similar structure as the ``mate in two'' puzzle:
we say a real function $f(x)$ has a limit at $x_0$ if there exists a real $a$ such that
for all $\varepsilon>0$ there exists $\delta>0$, provided 
$0<|x-x_0|<\delta$ implies $|f(x)-a|<\varepsilon$. Or, symbolically:
$$\exists a\in {\bf R}\hskip0.7cm
\forall \varepsilon >0\hskip0.7cm
\exists \delta >0\hskip0.7cm
\forall x: 0<|x-x_0|<\delta\hskip1cm
|f(x)-a|<\varepsilon.
$$

For instance, let us see how we can construct a concrete example similarly to the diagram for ``mate in two'':
Let us prove that the function $f(x)=x^2$ has a limit at $x=3$.

The limit is clearly $a=9$, this will be the first ``move'' for ``White''. Let us assume that ``Black'' responds
by using ``move'' $\varepsilon=1$. In this case the inequality $|x^2-9|<1$ implies the equivalent inequalities $8<x^2<10$ which
clearly hold when $|x-3|<3-\sqrt{8}$. That is, for this concrete epsilon the choice $\delta=3-\sqrt8$ is suitable.
For instance if $x=2.9$, then the difference of this result from 3 is less than
$3-\sqrt8\approx1.7$, and indeed, $2.9^2=8.41$, whose difference from $9$ is less than 1.
Or, going one step further, for $x=2.95$ the difference from 3 is even less: in this case $2.95^2=8.7025$ will
be even closer to 9. Obviously a number of other (infinite) instances can be mentioned.

Now let $\varepsilon=\frac{1}{10}$. In this case the inequality $|x^2-9|<\frac{1}{10}$ yields the equivalent inequalities $8.9<x^2<9.1$ 
which clearly hold if $|x-3|<3-\sqrt{8.9}$. That is, for this concrete epsilon the choice
$\delta=3-\sqrt{8.9}$ is suitable. Again, here we can do one or more concrete attempts: e.g. if $x=2.99$,
then its difference from 3 is less than $3-\sqrt{8.9}\approx0.016$, and indeed, $2.99^2=8.9401$, which
differs from $9$ in less than $\frac{1}{10}$.

It is clear that this kind of strategy can be played for any small positive epsilon. For an arbitrary ``black'' ``move''
$\varepsilon$ from the inequality $|x^2-9|<\varepsilon$ the equivalent inequalities $9-\varepsilon<x^2<9+\varepsilon$
hold which are necessarily true if $|x-3|<3-\sqrt{9-\varepsilon}$. That is, for any kind of attempt of ``Black's'',
``White's'' ``response'' $\delta=3-\sqrt{9-\varepsilon}$ will ensure the expected inequality described in the
definition, that is, the ``checkmate'' is guaranteed (see Fig.~\ref{fig:limit-function-graph}).
\begin{figure}
\centering
\begin{tabular}{ccccccccc}
$\exists a$ &                                       & $\forall\varepsilon$ &                         & $\exists\delta$ &            & $\forall x$ &               &    \\
\hline
     &                                          &                           &              &                    &  \rotatebox[origin=c]{5}{$\rightarrow$} & $2.9$ & \rotatebox[origin=c]{-20}{$\rightarrow$} &\\
     &  \rotatebox[origin=c]{10}{$\rightarrow$} & $1$           & $\rightarrow$ & $3-\sqrt8$ &  \rotatebox[origin=c]{0}{$\rightarrow$}&  $2.95$ & \rotatebox[origin=c]{-10}{$\rightarrow$} &\\
     &                                          &                           &               &                    &  $\vdots$                             & & &\\
$9$  &  \rotatebox[origin=c]{0}{$\rightarrow$}  & $\frac{1}{10}$ & $\rightarrow$ & $3-\sqrt{8.9}$ & \rotatebox[origin=c]{0}{$\rightarrow$}& $2.99$ & \rotatebox[origin=c]{0}{$\rightarrow$} & $|x^2-9|<\varepsilon$\\
     &  $\vdots$                                &                            &              &                        & $\vdots$ & & & \\
     &  \rotatebox[origin=c]{-10}{$\rightarrow$}& arb.& $\rightarrow$ & $3-\sqrt{9-\varepsilon}$  & $\rightarrow$ & $|x-3|<\delta$ & \rotatebox[origin=c]{20}{$\rightarrow$} &\\
\end{tabular}
\caption{Limit of a function, example graph}
\label{fig:limit-function-graph}
\end{figure}

\section{Remarks}
\begin{enumerate}
\item
Another idea could be that the structure $\exists\forall$ is demonstrated with ``mate in two''
and the structure $\exists\forall\exists\forall$ with ``mate in three''. In this case we can explain
this situation by the logical representation ``there exists a move for White such that for all
moves for Black, White \textit{may} deliver checkmate'', or by ``there exists a move for White such that for all
moves for Black, White has a move such that for all moves for Black, White \textit{may} deliver checkmate''.
In this approach it can be beneficial to not write such moves down which do not happen in real life
(i.e., the imaginary moves of the black king and the white rook). For those who are more familiar
with chess, this approach can be more convenient, but from the educational point of view it seems
to be simpler to do the demonstration with fewer moves.
\item
A widely used method to explain the classic definition of the limit is to \textit{deny} it.
What is more, the non-convergent sequences have their own names, notably,  \textit{divergent}
sequences. Another advantage of the analogy with the chess game that denial of the possibility
of checkmating is a natural situation. That is, an explanation of the impossibility of ``mate in two''
can be worded like ``for all moves for White, there exists a move for Black such that for all
moves for White, Black will not be checkmated''. Obviously this structure will be analogous to
the classic definition of the divergent sequences.
\item
One may think of defining the convergence by using the quantifiers with a more rigorous syntax. For example,
\cite[p.~48]{vajda} uses the notation
$$\displaystyle\mathop{\exists}_{\substack{a\in\bf R}}\ \mathop{\forall}_{\substack{\varepsilon\in\bf R\\ \varepsilon>0}}\ 
\displaystyle\mathop{\exists}_{\substack{M\in\bf R}}\ \mathop{\forall}_{\substack{x\in\bf R}}
\left(\left(x\geq M\right) \Rightarrow \left(\left|f(x)-a\right|<\varepsilon\right)\right)$$
to define the limit of a function at infinity to enable using Mathematica and the Theorema system
to prove certain properties strictly. Clearly this notation is beneficial
to obtain the negation of the definition as well, that is, since
$$A\Rightarrow B \equiv \neg A \vee B,$$
its negation is
$$\neg \left(A\Rightarrow B\right) \equiv \neg \left(\neg A \vee B\right) \equiv A \wedge \neg B,$$
hence the negation of the definition of the convergence will be
$$\displaystyle\mathop{\forall}_{\substack{a\in\bf R}}\ \mathop{\exists}_{\substack{\varepsilon\in\bf R\\ \varepsilon>0}}\ 
\displaystyle\mathop{\forall}_{\substack{M\in\bf R}}\ \mathop{\exists}_{\substack{x\in\bf R}}
\left(\left(x\geq M\right) \wedge \left(\left|f(x)-a\right|\geq\varepsilon\right)\right)$$
which could be rewritten as
$$\displaystyle\mathop{\forall}_{\substack{a\in\bf R}}\ \mathop{\exists}_{\substack{\varepsilon\in\bf R\\ \varepsilon>0}}\ 
\displaystyle\mathop{\forall}_{\substack{M\in\bf R}}\ \mathop{\exists}_{\substack{x\in\bf R\\ x\geq M}}
\left(\left|f(x)-a\right|\geq\varepsilon\right)$$
or---omitting some self-explanatory details---as
$$\displaystyle\mathop{\forall}_{\substack{a\in\bf R}}\ \mathop{\exists}_{\substack{\varepsilon>0}}\ 
\displaystyle\mathop{\forall}_{\substack{M\in\bf R}}\ \mathop{\exists}_{\substack{x\geq M}}
\left|f(x)-a\right|\geq\varepsilon.$$
On one hand this formula seems easy to interpret and memorize. On the other hand it seems difficult to
find easy methods to rewrite formulas mechanically, that is, manipulating on quantifier blocks
remains a challenge for most students.
\item
If the chess game appears to be difficult among some students, alternatively Bachet's game \cite{bachet} with 10 tokens on a table can
be considered, here the turns are taken by removing any number of tokens between 1 and 3 from the table. In Fig.~\ref{fig:bachet} Alice
removes 2 tokens in her first turn and she has a winning strategy to finally remove the last tokens as well.
\begin{figure}
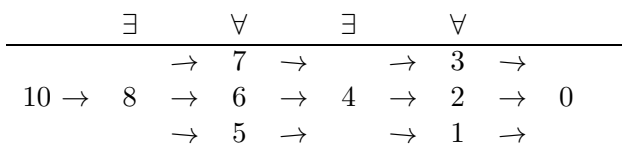

\centering
\begin{tabular}{ccccccccccc}
         & $\exists$ &                                         & $\forall$ &                                         & $\exists$ &                                       & $\forall$ &               &    \\
\hline
         &         & \rotatebox[origin=c]{5}{$\rightarrow$} & 7          & \rotatebox[origin=c]{-5}{$\rightarrow$} &           & \rotatebox[origin=c]{5}{$\rightarrow$} & 3          & \rotatebox[origin=c]{-5}{$\rightarrow$} &           &\\
10 $\to$ & 8       & \rotatebox[origin=c]{0}{$\rightarrow$} & 6          & \rotatebox[origin=c]{0}{$\rightarrow$}  & 4         & \rotatebox[origin=c]{0}{$\rightarrow$} & 2          & \rotatebox[origin=c]{0}{$\rightarrow$}  &  0        &\\
         &         & \rotatebox[origin=c]{-5}{$\rightarrow$}& 5          & \rotatebox[origin=c]{5}{$\rightarrow$}  &           & \rotatebox[origin=c]{-5}{$\rightarrow$} & 1         & \rotatebox[origin=c]{5}{$\rightarrow$} &           &\\
\end{tabular}
\caption{Winning strategy to control the amount of remaining tokens in Bachet's game}
\label{fig:bachet}
\end{figure}

\end{enumerate}
\section*{Acknowledgments}
The author thanks Stephen Jull, R\'obert Vajda and Peter Mayerhofer for comments that greatly improved the manuscript.

\bibliography{bibliography}{}
\bibliographystyle{apalike}

\end{document}